\documentclass{amsart}

\usepackage{amssymb,latexsym,amsfonts,amsmath}
\usepackage{multicol} 
\usepackage{graphicx}
\include{diagrams}
\usepackage{eso-pic}
\usepackage{epsfig}
\usepackage{graphicx}
\usepackage{color}
\usepackage{tikz}
\usetikzlibrary{automata}
\usetikzlibrary{shapes}
\usepackage{amssymb,latexsym,amsfonts,amsmath}
\usepackage{graphicx}
\newtheorem{theorem}{Theorem}[section]
\newtheorem{lemma}[theorem]{Lemma}

\newtheorem{definition}[theorem]{Definition}

\newtheorem{remark}[theorem]{Remark}
\numberwithin{equation}{section}

\newcommand{\R}{{\mathbb{R}}}

\newcommand{\X}{{\mathbf{X}}}

\begin{document}

\begin{abstract}
Stability is arguably one of the core concepts upon which our understanding of dynamical and control systems has been built. The related notion of incremental stability, however, has received much less attention until recently, when it was successfully used as a tool for the analysis and design of intrinsic observers, output regulation of nonlinear systems, frequency estimators, synchronization of coupled identical dynamical systems, symbolic models for nonlinear control systems, and bio-molecular systems. However, most of the existing controller design techniques provide controllers enforcing stability rather than incremental stability. Hence, there is a growing need to extend existing methods or develop new ones for the purpose of designing incrementally stabilizing controllers. In this paper, we develop a backstepping design approach for incremental stability. The effectiveness of the proposed method is illustrated by synthesizing a controller rendering a synchronous generator incrementally stable.
\end{abstract}

\title[Backstepping design for incremental stability]{Backstepping design for incremental stability}
\thanks{This work has been partially supported by the National Science Foundation award 0717188 and 0820061.}
\author[majid zamani]{Majid Zamani} 
\author[paulo tabuada]{Paulo Tabuada}
\address{Department of Electrical Engineering\\
University of California at Los Angeles,
Los Angeles, CA 90095}
\email{\{zamani, tabuada\}@ee.ucla.edu}
\urladdr{http://www.ee.ucla.edu/~zamani}
\urladdr{http://www.ee.ucla.edu/~tabuada}
\maketitle

\section{Introduction}
Stability is a property of dynamical systems comparing trajectories with an equilibrium point or with a particular trajectory. Incremental stability is a stronger property comparing arbitrary trajectories with themselves, rather than with an equilibrium point or with a particular trajectory. It is well-known that for linear systems incremental stability is equivalent to stability. For nonlinear systems, incremental stability is a stronger property requiring separate concepts and techniques for its study. 

%Incremental global asymptotic stability ($\delta$-GAS) and incremental input-to-state stability ($\delta$-ISS) are two common notions of incremental stability which have been successfully used by several researchers. 

The notion of incremental stability has a long history that can be traced back to the work of Zames in the 60's, as described in \cite{zames}. In \cite{zames1}, incremental stability is introduced and studied under the input-output setting where  control systems are regarded as operators mapping input signals to output signals. Incremental stability then arises naturally by considering the Lipschitz constant of the operator. A modern treatment of incremental stability, based on Lyapunov methods, appeared only recently in~\cite{angeli} where incremental global asymptotic stability ($\delta$-GAS) and incremental input-to-state stability ($\delta$-ISS) were defined in a state-space setting. The notion of $\delta$-GAS was defined by requiring the Euclidean distance between two arbitrary system trajectories to converge\footnote{A suitable ``small overshoot'' requirement is also included in the definition.} to zero. Lyapunov characterizations of $\delta$-GAS and $\delta$-ISS were also given in \cite{angeli}. There are two other stability properties related to incremental stability that have an equally long, if not longer, history. 

The first is the notion of convergent system which, according to~\cite{pavlov1}, was introduced in the 60's by B. P. Demidovich in~\cite{demidovich1,demidovich}. A system is convergent if all the system trajectories converge$^1$ to a bounded trajectory. Furthermore, Demidovich also introduced a sufficient condition for a system to be convergent, called the Demidovich's condition in \cite{pavlov}. Since incremental stability requires every trajectory to converge to every other trajectory, an incrementally stable system is also a convergent system whenever a bounded trajectory exists.

The second stability property is contractivity and was introduced in the control community by Lohmiller and Slotine in~\cite{lohmiller} although it had been studied before in the mathematical community~\cite{jerome}. Rather than comparing trajectories, the notion of contracting system is infinitesimal and requires the decrease of a suitable quantity, defined through a Riemannian metric, along trajectories. The definition of contracting system can be seen as a generalization of the Demidovich's condition. Note that while the Demidovich's condition was introduced as a sufficient condition for a system to be convergent, its generalization in~\cite{lohmiller} was directly used as the definition of contracting system. 

Although both contractivity as well as convergence are coordinate independent properties, this is no longer the case with incremental stability. In this paper, however, we work with a variation of incremental stability that is coordinate invariant. This is achieved by no longer insisting on the distance between trajectories being measured by the Euclidean metric.

%Recent work related to contraction analysis can be found in \cite{aghannan, erin, pham, majid, sontag2}. 

%As showed in \cite{pavlov}, the property of convergence is coordinate independent, but $\delta$-GAS and $\delta$-ISS properties are coordinate dependent. On the other hand, $\delta$-GAS and $\delta$-ISS properties provide convergence of trajectories with respect to themselves, rather than only a bounded trajectory, the case in convergence property. The contraction property, as we show in this paper, is coordinate independent, but it implies convergence of trajectories with respect to themselves and with respect to some metric, not necessarily Euclideam metric, the case in convergence, $\delta$-GAS and $\delta$-ISS properties.

The number of applications of incremental stability has increased in the past years. Examples include intrinsic observer design \cite{aghannan}, consensus problems in complex networks \cite{wang}, output regulation of nonlinear systems \cite{pavlov}, design of frequency estimators \cite{sharma1}, synchronization of coupled identical dynamical systems \cite{russo1}, construction of symbolic models for nonlinear control systems \cite{pola, pola1, girard2}, and the analysis of bio-molecular systems \cite{russo}. Our motivation comes from symbolic control where incremental stability was identified as a key property enabling the construction of finite abstractions of nonlinear control systems~\cite{pola, pola1, girard2,majid}. Hence, there is a growing need for design methods providing controllers enforcing incremental stability since most of the existing design methods guarantee stability rather than incremental stability. 

Related work includes the recent design results enforcing the convergent system property through the solution of linear matrix inequalities (LMIs) \cite{pavlov,pavlov2,wouw}. In contrast, the results presented in this paper do not require the solution of LMIs and existence of controllers is always guaranteed. Backstepping design methods for incremental global asymptotic stability\footnote{Understood in the sense of Definition~\ref{delta_GAS}.} for parametric-strict-feedback\footnote{See equation~(\ref{system1}) or~\cite{miroslav} for a definition.} systems were proposed before in \cite{jouffroy1, sharma2}, and \cite{sharma}. In this paper, we generalize the results in \cite{jouffroy1, sharma2}, and \cite{sharma} by:
\begin{itemize}
\item[1)] developing a backstepping design method providing controllers enforcing incremental input-to-state stability\footnote{Understood in the sense of Definition~\ref{delta_ISS}.} and not simply incremental global asymptotic stability;
\item[2)] enlarging the class of control systems from parametric-strict-feedback to strict-feedback\footnote{See equation~(\ref{system6}) or~\cite{miroslav} for a definition.} form.  
\end{itemize}
%develop a backstepping design method providing controllers enforcing $\delta$-GAS and $\delta$-ISS for more general class of control systems, called strict-feedback systems \cite{miroslav}. 
The proposed approach was inspired by the original backstepping method described, for example, in \cite{miroslav}. Like the original backstepping method, which provides a recursive way of constructing controllers as well as Lyapunov functions, the approach proposed in this paper provides a recursive way of constructing controllers as well as contraction metrics. Our design approach is illustrated by designing a controller rendering a synchronous generator incrementally stable.

\section{Control Systems and Stability Notions}
\subsection{Notation} 
The symbols $\mathbb{R}$, $\mathbb{R}^+$ and $\mathbb{R}_0^+$ denote the set of real, positive, and nonnegative real numbers, respectively. The symbols $I_m$, and $0_m$ denote the identity and zero matrices in $\R^{m\times{m}}$. Given a vector $x\in\mathbb{R}^{n}$, we denote by $x_{i}$ the \mbox{$i$--th} element of $x$, and by $\Vert x\Vert$ the Euclidean norm of $x$; we recall that \mbox{$\Vert x\Vert=\sqrt{x_1^2+x_2^2+...+x_n^2}$}. 
%A function \mbox{$f:[a,b]\rightarrow\mathbb{R}^n$} is said to be absolutely continuous on $[a,b]$ if for any $\varepsilon\in\mathbb{R^+}$ there exists $\delta\in\mathbb{R}^+$ so that for every $k\in\mathbb{N}$ and for every sequence of points \mbox{$a\leq{a_1}<b_1<a_2<b_2<\ldots<a_k<b_k\leq{b}$}, if \mbox{$\sum_{i=1}^m(b_i-a_i)<\delta$} then \mbox{$\sum_{i=1}^m\Vert f(b_i)-f(a_i)\Vert<\varepsilon$}. A function \mbox{$f:]a,b[\rightarrow\mathbb{R}^n$} is said to be locally absolutely continuous if the restriction of $f$ to any compact subset of $]a,b[$ is absolutely continuous. 
Given a measurable function \mbox{$f:\mathbb{R}_{0}^{+}\rightarrow\mathbb{R}^n$}, the (essential) supremum of $f$ is denoted by $\Vert f\Vert_{\infty}$; we recall that \mbox{$\Vert f\Vert_{\infty}:=\text{(ess)sup}\{\Vert f(t)\Vert,t\geq0\}$}. 
%$f$ is essentially bounded if $\Vert{f}\Vert_{\infty}<\infty$. For a given time $\tau\in\mathbb{R}^+$, define $f_{\tau}$ so that $f_{\tau}(t)=f(t)$, for any $t\in[0,\tau)$, and $f(t)=0$ elsewhere; $f$ is said to be locally essentially bounded if for any $\tau\in\mathbb{R}^+$, $f_{\tau}$ is essentially bounded. 
A continuous function \mbox{$\gamma:\mathbb{R}_{0}^{+}\rightarrow\mathbb{R}_{0}^{+}$}, is said to belong to class $\mathcal{K}$ if it is strictly increasing and \mbox{$\gamma(0)=0$}; $\gamma$ is said to belong to class $\mathcal{K}_{\infty}$ if \mbox{$\gamma\in\mathcal{K}$} and $\gamma(r)\rightarrow\infty$ as $r\rightarrow\infty$. A continuous function \mbox{$\beta:\mathbb{R}_{0}^{+}\times\mathbb{R}_{0}^{+}\rightarrow\mathbb{R}_{0}^{+}$} is said to belong to class $\mathcal{KL}$ if, for each fixed $s$, the map $\beta(r,s)$ belongs to class $\mathcal{K}_{\infty}$ with respect to $r$ and, for each fixed nonzero $r$, the map $\beta(r,s)$ is decreasing with respect to $s$ and $\beta(r,s)\rightarrow0$ as \mbox{$s\rightarrow\infty$}. If $\phi:\R^n\to \R^n$ is a smooth function with a smooth inverse, called a diffeomorphism, and if $X:\R^n\to \R^n$ is a smooth map, we denote by $\phi_*X$ the map defined by $(\phi_*X)(y)=\frac{\partial \phi}{\partial x}\big\vert_{x=\phi^{-1}(y)}X\circ\phi^{-1}(y)$. Let now $G:\R^n\to \R^{n\times n}$ be a smooth map. The notation $\phi^*G:\R^n\to \R^{n\times n}$ denotes the smooth map $(\phi^*G)(x)=(\frac{\partial \phi}{\partial x})^TG(\phi(x))(\frac{\partial \phi}{\partial x})$. A Riemannian metric \mbox{$G:\mathbb{R}^n\rightarrow\mathbb{R}^{n\times n}$} is a smooth map on $\mathbb{R}^n$ such that, for any \mbox{$x\in\mathbb{R}^n$}, $G(x)$ is a symmetric positive definite matrix \cite{lee}. For any \mbox{$x\in\mathbb{R}^n$} and smooth functions \mbox{$I,J:\mathbb{R}^n\rightarrow\mathbb{R}^n$}, one can define the scalar function \mbox{$\langle{I},J\rangle_G$} as $I^T(x)G(x)J(x)$. We will still use the notation $\langle{I},J\rangle_{G}$ to denote $I^TGJ$ even if $G$ does not represent any Riemannian metric. A function \mbox{$\mathbf{d}:\R^n\times \R^n\rightarrow\mathbb{R}_{0}^{+}$} is a metric on $\R^n$ if for any $x,y,z\in\R^n$, the following three conditions are satisfied: i) $\mathbf{d}(x,y)=0$ if and only if $x=y$; ii) $\mathbf{d}(x,y)=\mathbf{d}(y,x)$; and iii) $\mathbf{d}(x,z)\leq\mathbf{d}(x,y)+\mathbf{d}(y,z)$. We use the pair $\left(\R^n,\mathbf{d}\right)$ to denote a metric space $\R^n$ equipped with the metric $\mathbf{d}$. We use the notation $\mathbf{d}_G$ to denote the Riemannian distance function provided by the Riemannian metric $G$ \cite{lee}.

\subsection{Control Systems\label{II.B}}

The class of control systems that we consider in this paper is formalized in
the following definition.
\begin{definition}
\label{Def_control_sys}A \textit{control system} is a quadruple:
\[
\Sigma=(\mathbb{R}^{n},\mathsf{U},\mathcal{U},f),
\]
where:
\begin{itemize}
\item $\mathbb{R}^{n}$ is the state space;
\item $\mathsf{U}\subseteq\mathbb{R}^{m}$ is the input set which is compact and convex;
\item $\mathcal{U}$ is the set of all measurable functions of time from intervals of the form \mbox{$]a,b[\subseteq\mathbb{R}$} to $\mathsf{U}$ with $a<0$ and $b>0$; 
\item \mbox{$f:\mathbb{R}^{n}\times \mathsf{U}\rightarrow\mathbb{R}^{n}$} is a continuous map
satisfying the following Lipschitz assumption: for every compact set
\mbox{$Q\subset\mathbb{R}^{n}$}, there exists a constant $Z\in\mathbb{R}^+$ such that $\Vert
f(x,u)-f(y,u)\Vert\leq Z\Vert x-y\Vert$ for all $x,y\in Q$ and all $u\in \mathsf{U}$.
\end{itemize}
\end{definition}

A curve \mbox{$\xi:]a,b[\rightarrow\mathbb{R}^{n}$} is said to be a
\textit{trajectory} of $\Sigma$ if there exists $\upsilon\in\mathcal{U}$
satisfying:
\begin{equation}
\dot{\xi}(t)=f\left(\xi(t),\upsilon(t)\right),\label{eq0}
\end{equation}
for almost all $t\in$ $]a,b[$. 
%Although we have defined trajectories over open domains, we shall refer to
%trajectories \mbox{${\xi:}[0,\tau]\rightarrow\mathbb{R}^{n}$} defined on closed
%domains $[0,\tau],$ $\tau\in\mathbb{R}^{+}$ with the understanding of the
%existence of a trajectory \mbox{${\xi}^{\prime}:]a,b[\rightarrow\mathbb{R}^{n}$}
%such that \mbox{${\xi}={\xi}^{\prime}|_{[0,\tau]}$}. 
We also write $\xi_{x\upsilon}(t)$ to denote the point reached at time $t$
under the input $\upsilon$ from initial condition $x=\xi_{x\upsilon}(0)$; this point is
uniquely determined, since the assumptions on $f$ ensure existence and
uniqueness of trajectories \cite{sontag1}. 
%We also denote an autonomous system $\Sigma$ with no control inputs by $\Sigma=(\mathbb{R}^n,f)$. 
A control system $\Sigma$ is said to be forward complete if every trajectory is defined on an interval of the form $]a,\infty[$. Sufficient and necessary conditions for a system to be forward complete can be found in \cite{sontag}. A control system $\Sigma$ is said to be smooth if $f$ is an infinitely differentiable function of its arguments.

\subsection{Stability notions}
We start by introducing the following definitions which were inspired by the notions of incremental global asymptotic stability ($\delta$-GAS) and incremental input-to-state stability ($\delta$-ISS) presented in \cite{angeli}.
\begin{definition}
A control system $\Sigma$ is incrementally globally asymptotically stable ($\delta_\exists$-GAS) if it is forward complete and there exist a metric $\mathbf{d}$ and a $\mathcal{KL}$ function $\beta$ such that for any $t\in{\mathbb{R}_0^+}$, any $x,x'\in{\mathbb{R}^n}$ and any $\upsilon\in\mathcal{U}$ the following condition is satisfied:
\begin{equation}
\mathbf{d}\left(\xi_{x\upsilon}(t),\xi_{x'\upsilon}(t)\right) \leq\beta\left(\mathbf{d}\left(x,x'\right),t\right). \label{delta_GAS}%
\end{equation}
\end{definition}
While $\delta$-GAS, as defined in~\cite{angeli}, requires the metric $\mathbf{d}$ to be the Euclidean metric, Definition~\ref{delta_GAS} only requires the existence of a metric; hence, the existential quantifier in the acronym $\delta_\exists$-GAS. We note that while $\delta$-GAS is not invariant under changes of coordinates, $\delta_\exists$-GAS is. If $\phi:\R^n\to \R^n$ is a bijective change of coordinates, inequality~(\ref{delta_GAS}) transforms under $\phi$ to: 
$$\mathbf{d}'\left(\phi\circ\xi_{x\upsilon}(t),\phi\circ\xi_{x'\upsilon}(t)\right) \leq\beta\left(\mathbf{d}'\left(\phi(x),\phi(x')\right),t\right),$$
where $\mathbf{d}'(y,y')=\mathbf{d}(\phi^{-1}(y),\phi^{-1}(y'))$. We shall return to the comparison between $\delta_\exists$-GAS and $\delta$-GAS in Subsection~\ref{SSec:Descriptions}. Nevertheless, when the origin is an equilibrium point for $\Sigma$ and the map $\psi:\R^n\to \R_0^+$, defined by $\psi(x)=\mathbf{d}(x,0)$, is continuous and radially unbounded\footnote{Under the stated assumptions it can be shown that $\underline{\alpha}(\Vert x\Vert)\le \psi(x)\le\overline{\alpha}(\Vert x\Vert)$ for $\mathcal{K}_\infty$ functions $\underline{\alpha}$ and $\overline{\alpha}$.}, both $\delta_\exists$-GAS and $\delta$-GAS imply global asymptotic stability.

\begin{definition}
\label{dISS}
A control system $\Sigma$ is incrementally input-to-state stable ($\delta_\exists$-ISS) if it is forward complete and there exist a metric $\mathbf{d}$, a $\mathcal{KL}$ function $\beta$, and a $\mathcal{K}_{\infty}$ function $\gamma$ such that for any $t\in{\mathbb{R}_0^+}$, any $x,x'\in{\mathbb{R}^n}$, and any $\upsilon$, ${\upsilon}'\in\mathcal{U}$ the following condition is satisfied:
\begin{equation}
\mathbf{d}\left(\xi_{x\upsilon}(t),\xi_{x'{\upsilon}'}(t)\right) \leq\beta\left(\mathbf{d}\left(x,x'\right),t\right)+\gamma\left(\left\Vert{\upsilon}-{\upsilon}'\right\Vert_{\infty}\right). \label{delta_ISS}%
\end{equation}
\end{definition}

By observing (\ref{delta_GAS}) and (\ref{delta_ISS}), it is readily seen that $\delta_\exists$-ISS implies $\delta_\exists$-GAS while the converse is not true in general. Moreover, whenever the metric $\mathbf{d}$ is the Euclidean metric, $\delta_\exists$-ISS becomes $\delta$-ISS as defined in~\cite{angeli}. We note that while $\delta$-ISS is not invariant under changes of coordinates, $\delta_\exists$-ISS is. Once again, although $\delta_\exists$-ISS is not equivalent to $\delta$-ISS, both notions imply input-to-state stability whenever the origin is an equilibrium point for $\Sigma$ and the map $\psi:\R^n\to \R_0^+$, defined by $\psi(x)=\mathbf{d}(x,0)$, is continuous and radially unbounded.

\subsection{Descriptions of incremental stability}
\label{SSec:Descriptions}

One of the methods for checking incremental stability properties consists in using Lyapunov functions. The Lyapunov characterizations of $\delta$-GAS and $\delta$-ISS were developed in \cite{angeli}. In this paper we follow an alternative approach based on contraction metrics describing $\delta_\exists$-GAS and $\delta_\exists$-ISS properties, when $\mathbf{d}$ is the Riemannian distance function. The notion of contraction metric was popularized in control theory by the work of Lohmiller and Slotine \cite{lohmiller} that relies on variational systems.

The variational system associated with a smooth control system $\Sigma=(\mathbb{R}^n,\mathsf{U},\mathcal{U},f)$, when we only have state variations, is given by the differential equation: 
\begin{equation}
\label{variational1}
\frac{d}{dt}(\delta{\xi})=\frac{\partial{f}}{\partial{x}}\bigg |_{x=\xi\atop u=\upsilon}\delta{\xi},
\end{equation}
for any $\upsilon\in\mathcal{U}$ and where $\delta{\xi}$ is the variation\footnote{The variation $\delta{\xi}$ can be formally defined by considering a family of trajectories $\xi_{x\upsilon}(t,\epsilon)$ parametrized by $\epsilon\in\mathbb{R}$. The variation of the state is then $\delta{\xi}=\frac{\partial{\xi}_{x\upsilon}}{\partial{\epsilon}}$.} of a trajectory $\xi$ of $\Sigma$. More details about the variational system can be found in \cite{crouch}. Similarly, the variational system associated with a smooth control system $\Sigma=(\mathbb{R}^n,\mathsf{U},\mathcal{U},f)$, when we have state and input variations, is given by the differential equation:
\begin{equation}
\label{variational2}
\frac{d}{dt}(\delta{\xi})=\frac{\partial{f}}{\partial{x}}\bigg |_{x=\xi\atop u=\upsilon}\delta{\xi}+\frac{\partial{f}}{\partial{u}}\bigg |_{x=\xi\atop u=\upsilon}\delta{\upsilon},
\end{equation}
where $\delta{\xi}$ and $\delta{\upsilon}$ are variations of a state trajectory $\xi$ and an input curve $\upsilon$ of $\Sigma$, respectively. 

The following definition is adapted from~\cite{lohmiller}:

\begin{definition}
\label{contraction}
Let $\Sigma=(\mathbb{R}^n,\mathsf{U},\mathcal{U},f)$ be a smooth control system equipped with a Riemannian metric $G$. The metric $G$ is said to be a contraction metric, with respect to states, for system $\Sigma$ if there exists some $\lambda\in\mathbb{R}^+$ such that:
\begin{equation}
\langle X,X \rangle_{F}\leq-\lambda\langle X, X\rangle_{G} \label{contractionGAS1}
\end{equation}
for $F(x,u)=\left(\frac{\partial{f}}{\partial{x}}\right)^TG(x)+G(x)\frac{\partial{f}}{\partial{x}}+\frac{\partial{G}}{\partial{x}}f(x,u)$, any $u\in\mathsf{U}$, $X\in\mathbb{R}^n$ and $x\in \R^n$, or equivalently:
\begin{eqnarray}
X^T\left(\left(\frac{\partial{f}}{\partial{x}}\right)^TG(x)+G(x)\frac{\partial{f}}{\partial{x}}+\frac{\partial{G}}{\partial{x}}f(x,u)\right)X\leq-\lambda X^TG(x)X, \label{contractionGAS}
\end{eqnarray}
where the constant $\lambda$ is called contraction rate. 
\end{definition}
When the metric $G$ is constant, the condition (\ref{contractionGAS1}) or (\ref{contractionGAS}) is known as the Demidovich's condition \cite{pavlov}. It is shown in~\cite{pavlov} that such condition implies incremental stability and the convergent system property.

Note that the inequality (\ref{contractionGAS1}) or (\ref{contractionGAS}) implies: 
\begin{equation}
\label{dervar1}
\frac{d}{dt}\langle\delta{\xi},\delta{\xi}\rangle_{G}\leq-\lambda\langle\delta{\xi},\delta{\xi}\rangle_{G}, 
\end{equation}
when we only have state variations and $\delta\xi$ is the variation of a state trajectory $\xi$ of $\Sigma$. 

%We say that a smooth control system $\Sigma$ is an exponential state contraction if there exists a Riemannian metric $G$ such that $\Sigma$ is an exponential state contraction with respect to the metric $G$. Moreover, we call $G$ a state contraction metric.

The following theorem shows that the inequality (\ref{contractionGAS1}) implies $\delta_\exists$-GAS.
\begin{theorem}
\label{theoremGAS}
Let \mbox{$\Sigma=(\mathbb{R}^n,\mathsf{U},\mathcal{U},f)$} be a smooth control system equipped with a Riemannian metric $G$. If $G$ is a contraction metric, with respect to states, for $\Sigma$ and if $\left(\R^n,\mathbf{d}_G\right)$ is a complete metric space\footnote{A metric space $(\R^n,\mathbf{d})$ is said to be complete if every Cauchy sequence of points in $\R^n$ has a limit that is also in $\R^n$.}, then $\Sigma$ is $\delta_\exists$-GAS.
\end{theorem}

Different variations of this result appeared in \cite{lohmiller} and \cite{aghannan}; see \cite{aghannan} for a concise proof and~\cite{majid} for a more detailed treatment including the completeness assumption. It is also shown in~\cite{majid} that the additional assumption $\underline{\omega}\langle X,X\rangle_{I_n}\le \langle X,X\rangle_G\le \overline{\omega}\langle X,X\rangle_{I_n}$ for $\underline{\omega},\overline{\omega}\in \R^+$ leads to the stronger conclusion that $\Sigma$ is in fact $\delta$-GAS.

If, in addition to state variations, we also allow for input variations we obtain the notion of contraction metric with respect to states and inputs.
\begin{definition}
\label{contraction2}
Let $\Sigma=(\mathbb{R}^n,\mathsf{U},\mathcal{U},f)$ be a smooth control system on $\mathbb{R}^n$ equipped with a Riemannian metric $G$. The metric $G$ is said to be a contraction metric, with respect to states and inputs, for system $\Sigma$ if there exists some $\lambda\in\mathbb{R}^+$ and $\alpha\in\mathbb{R}_0^+$ such that:
\begin{eqnarray}
\label{contraction3ISS}
\langle{X},X\rangle_{F}+2\bigg\langle\frac{\partial{f}}{\partial{u}}Y,X\bigg\rangle_{G}\leq-\lambda\langle{X},X\rangle_{G}+\alpha{\langle{X},X\rangle^{\frac{1}{2}}_G}\langle{Y,Y}\rangle^{\frac{1}{2}}_{I_m}
\end{eqnarray}
for $F(x,u)=\left(\frac{\partial{f}}{\partial{x}}\right)^TG(x)+G(x)\frac{\partial{f}}{\partial{x}}+\frac{\partial{G}}{\partial{x}}f(x,u)$, any $X\in\mathbb{R}^n$, $x\in\R^n$, $u\in\mathsf{U}$, and $Y\in\R^m$, or equivalently:
\begin{eqnarray}
\nonumber
&X^T\left(\left(\frac{\partial{f}}{\partial{x}}\right)^TG(x)+G(x)\frac{\partial{f}}{\partial{x}}+\frac{\partial{G}}{\partial{x}}f(x,u)\right)X+2Y^T\left(\frac{\partial{f}}{\partial{u}}\right)^TG(x)X\\
&\leq-\lambda{X}^TG(x)X+\alpha(X^TG(x)X)^\frac{1}{2}(Y^TY)^\frac{1}{2}, \label{contraction2ISS}
\end{eqnarray} 
where the constant $\lambda$ is called contraction rate. 
\end{definition}
Note that the inequality (\ref{contraction3ISS}) or (\ref{contraction2ISS}) implies: 
\begin{equation}
\label{dervar2}
\frac{d}{dt}\langle\delta{\xi},\delta{\xi}\rangle_{G}\leq-\lambda\langle\delta{\xi},\delta{\xi}\rangle_{G}+\alpha{\langle{\delta\xi},\delta\xi\rangle^{\frac{1}{2}}_G}\langle{\delta\upsilon,\delta\upsilon}\rangle^{\frac{1}{2}}_{I_m}, 
\end{equation}
when we have state and input variations and where $\delta\xi$ and $\delta\upsilon$ are variations of a state trajectory $\xi$ and an input curve $\upsilon$ of $\Sigma$. 

%We say that a smooth control system $\Sigma$ is an exponential state-input contraction if there exists a Riemannian metric $G$ such that $\Sigma$ is an exponential state-input contraction with respect to the metric $G$. Moreover, we call $G$ a state-input contraction metric.

The following theorem shows that the inequality (\ref{contraction3ISS}) implies $\delta_\exists$-ISS.
\begin{theorem}[\cite{majid}]
\label{theoremISS}
Let $\Sigma=(\mathbb{R}^n,\mathsf{U},\mathcal{U},f)$ be a smooth control system on $\mathbb{R}^n$ equipped with a Riemannian metric $G$. If the metric $G$ is a contraction metric, with respect to states and inputs, for system $\Sigma$ and $\left(\R^n,\mathbf{d}_G\right)$ is a complete metric space, then $\Sigma$ is $\delta_\exists$-ISS.
\end{theorem}

Similarly to contraction metrics with respect to states, it shown in~\cite{majid} that $\Sigma$ is $\delta$-ISS when the inequalities $\underline{\omega}\langle X,X\rangle_{I_n}\le \langle X,X\rangle_G\le \overline{\omega}\langle X,X\rangle_{I_n}$ are satisfied for $\underline{\omega},\overline{\omega}\in \R^+$.

In the next section, we propose a backstepping design procedure to render control systems incrementally stable. 

\section{Backstepping Design Procedure}\label{back}
The method described here was inspired by the original backstepping described, for example, in \cite{miroslav}. Consider the class of control systems $\Sigma=(\R^n,\mathsf{U},\mathcal{U},f)$ with $f$ of the parametric-strict-feedback form \cite{miroslav}:
\begin{eqnarray}
                \begin{array}{ccl}
                 f_1(x,u)&=&h_1(x_1)+b_1x_2,\\
                 f_2(x,u)&=&h_2(x_1,x_2)+b_2x_3,\\
                 &\vdots&\\
                 f_{n-1}(x,u)&=&h_{n-1}(x_1,\cdots,x_{n-1})+b_{n-1}x_n,\\
                 f_{n}(x,u)&=&h_{n}(x)+g(x)u,\\
                \end{array}
\label{system1}
\end{eqnarray}
where $x\in\R^n$ is the state and $u\in\mathsf{U}\subset\R$ is the control input. The functions \mbox{$h_i:\mathbb{R}^{i}\rightarrow\mathbb{R}$}, for \mbox{$i=1,\ldots,n$}, and \mbox{$g:\mathbb{R}^{n}\rightarrow\mathbb{R}$} are smooth, \mbox{$g(x)\neq{0}$} over the domain of interest, and $b_i\in\mathbb{R}$, for \mbox{$i=1,\ldots,n$}, are nonzero constants. 
We can now state one of the main results, describing a backstepping controller for control system (\ref{system1}).
\begin{theorem}
\label{theorem1}
For any control system $\Sigma=(\R^n,\mathsf{U},\mathcal{U},f)$ with $f$ of the form~(\ref{system1}) and for any $\lambda\in\mathbb{R}^+$, the state feedback control law:
\begin{eqnarray}
k(x,\widehat{u})&=&\frac{1}{g(x)}\bigg[k_{n}(x)-h_n(x)\bigg]+\frac{1}{g(x)}\widehat{u}, \label{inputor}
\end{eqnarray}
where 
\begin{eqnarray}
\nonumber
k_{l}(x,\widehat{u})&=&-b_{l-1}\left(x_{l-1}-\phi_{l-2}(x)\right)-\frac{\lambda}{2}\left(x_{l}-\phi_{l-1}(x)\right)+\frac{\partial{\phi}_{l-1}}{\partial{x}}f(x,k(x,\widehat{u})),~\text{for}~l=1,\cdots,n,\\ \notag
\phi_{l}(x)&=&\frac{1}{b_{l}}\bigg[k_{l}(x)-h_{l}(x)\bigg],~\text{for}~l=1,\cdots,n-1,\\ \notag && \phi_{-1}(x)=\phi_0(x)=0~\forall x\in\R^n,~b_0=0,~\text{and}~ x_{0}=0,  \label{input}
\end{eqnarray}
renders the control system $\Sigma$ $\delta_\exists$-GAS.
\end{theorem}

\begin{proof}
Consider the following system:
%\begin{eqnarray}
%    \Sigma_l:\left\{
%                \begin{array}{ccl}
%                 \dot{\xi}_1&=&f_1(\xi_1)+b_1\xi_2,\\
%                 &\vdots&\\
%                 \dot{\xi}_{l-1}&=&f_{l-1}(\xi_1,\cdots,\xi_{l-1})+b_{l-1}\xi_{l},\\
%                 \dot{\xi}_{l}&=&k_l(\xi_1,\cdots,\xi_{l}),\\
%                \end{array}
%                \right.
%\label{system2}
%\end{eqnarray}
\begin{equation}
    \Sigma_l:\left\{
                \begin{array}{ll}
                 \dot{\eta}_l=F_l(\eta_l)+B_l\xi_l,\\
                 \dot{\xi}_l=k_l(\eta_l,\xi_l),\\
                \end{array}
                \right.
\label{system3}
\end{equation}
%or equivalently: \r{***Decide which one you want to use. Don't start a proof defining a control system in two different ways!***}
%\begin{equation}
%    \Sigma_l:\left\{
%                \left[\begin{array}{c}
%                 \dot{\eta}_l\\
%                 \dot{\xi}_l\\
%                \end{array}
%                \right]=F_{l+1}(\eta_{l+1})+B_{l+1}\phi_{l},
%                \right.
%\label{system31}
%\end{equation}
where \mbox{$\eta_l=[\xi_1,\cdots,\xi_{l-1}]^T$}, $B_l=[0,\cdots,0,b_{l-1}]^T\in\mathbb{R}^{l-1}$, \mbox{$z_l=[y_l^T~x_l]^T\in\mathbb{R}^l$} is the state of $\Sigma_l$, $y_l=[x_1,\cdots,x_{l-1}]^T$, and \mbox{$F_l(y_l)=[f_1(x,u),\cdots,f_{l-2}(x,u),h_{l-1}(x_1,\cdots,x_{l-1})]^T$}. By using induction on $l$, we show that the metric $G_l$, defined by: 
%\begin{small}
\begin{equation}
G_l(y_l)=\left[ {\begin{array}{cc}
 G_{l-1}(y_{l-1})+\left(\frac{\partial{\phi_{l-1}}}{\partial{y_l}}\right)^T\frac{\partial{\phi_{l-1}}}{\partial{y_l}}&-\left(\frac{\partial{\phi_{l-1}}}{\partial{y_l}}\right)^T  \\
 -\frac{\partial{\phi_{l-1}}}{\partial{y_l}}&1  \\
 \end{array} } \right], \label{contraction_metric}
\end{equation}
%\end{small}
is a contraction metric, with respect to states, for the system (\ref{system3}) with contraction rate $\lambda$. For $l=1$, it can be easily checked that $G_1(y_1)=1$ is a contraction metric, with respect to states, with the contraction rate $\lambda$ for the scalar system: 
\begin{eqnarray}
\nonumber
\Sigma_1:\dot{\xi}_1=k_1(\xi_1)=-\frac{\lambda}{2}\xi_1. \nonumber
\end{eqnarray}
Assume that the metric $G_{k-1}$ is a contraction metric, with respect to states, for the system $\Sigma_{k-1}$, for some $2\leq k\leq n$, and with contraction rate $\lambda$. This implies:
\begin{footnotesize}
\begin{eqnarray}
\nonumber
\left[Y^T~X\right]\left(\left(\frac{\partial(F_k+B_k\phi_{k-1})}{\partial{y_k}}\right)^TG_{k-1}(y_{k-1})+G_{k-1}(y_{k-1})\frac{\partial(F_k+B_k\phi_{k-1})}{\partial{y_k}}+\frac{\partial{G}_{k-1}}{\partial{y_k}}\left(F_k+B_k\phi_{k-1}\right)\right)\left[{\begin{array}{c}Y\\X\\\end{array}}\right]\\ \label{contsub}
\leq-\lambda\left[Y^T~X\right]G_{k-1}(y_{k-1})\left[{\begin{array}{c}Y\\X\\\end{array}}\right], 
\end{eqnarray}
\end{footnotesize}
for any $Y\in\R^{k-2}$, and $X\in\R$.
Since the metric $G_{k-1}$ is only a function of $y_{k-1}=\left[x_1,\cdots,x_{k-2}\right]^T$, and the vector $B_k$ has zero entries except for the last entry, it can be easily shown that $\frac{\partial{G}_{k-1}}{\partial{y_k}}B_k=0_{k-1}$, and the inequality (\ref{contsub}) reduces to:
\begin{small}
\begin{eqnarray}
\nonumber
\left[Y^T~X\right]\left(\left(\frac{\partial(F_k+B_k\phi_{k-1})}{\partial{y_k}}\right)^TG_{k-1}(y_{k-1})+G_{k-1}(y_{k-1})\frac{\partial(F_k+B_k\phi_{k-1})}{\partial{y_k}}+\frac{\partial{G}_{k-1}}{\partial{y_k}}F_k\right)\left[{\begin{array}{c}Y\\X\\\end{array}}\right]\\ \label{contsub1}
\leq-\lambda\left[Y^T~X\right]G_{k-1}(y_{k-1})\left[{\begin{array}{c}Y\\X\\\end{array}}\right]. 
\end{eqnarray}
\end{small}
Now, we show that:
%\begin{small}
\begin{equation}
G_k(y_{k})=\left[ {\begin{array}{cc}
 G_{k-1}(y_{k-1})+\left(\frac{\partial{\phi_{k-1}}}{\partial{y_{k}}}\right)^T\frac{\partial{\phi_{k-1}}}{\partial{y_{k}}}&-\left(\frac{\partial{\phi_{k-1}}}{\partial{y_{k}}}\right)^T  \\
 -\frac{\partial{\phi_{k-1}}}{\partial{y_{k}}}&1  \\
 \end{array} } \right], \label{GASmetric}
\end{equation}
%\end{small}
is a contraction metric, with respect to states, for the system $\Sigma_k$.
%\begin{equation}
%    \Sigma_k:\left\{
%                \begin{array}{ll}
%                 \dot{\eta}=f(\eta)+b\xi_k,\\
%                 \dot{\xi}_k=k_k(\eta,\xi_k),\\
%                \end{array}
%                \right.
%\label{system4}
%\end{equation}
%where \mbox{$\eta=[\xi_1,\cdots,\xi_{k-1}]^T$}, $b=[0,\cdots,0,b_{k-1}]$, \mbox{$f(\eta)=[f_1(\xi_1)+b_1\xi_2,\cdots,f_{k-1}(\xi_1,\cdots,\xi_{k-1})]^T$}, and \mbox{$y=\eta_{y\upsilon}(0)$}. 
For any nonzero vector $\left[Y^T~X\right]^T\in\R^k$, we have:
%\begin{small}
\begin{eqnarray}
\nonumber
\left[Y^T~X\right]G_{k}(y_k)\left[{\begin{array}{c}Y\\X\\\end{array}}\right]&=&
\left[Y^T~X\right]\left[{\begin{array}{cc}
 G_{k-1}(y_{k-1})+\left(\frac{\partial{\phi_{k-1}}}{\partial{y_k}}\right)^T\frac{\partial{\phi_{k-1}}}{\partial{y_k}}&-\left(\frac{\partial{\phi_{k-1}}}{\partial{y_k}}\right)^T  \\
 -\frac{\partial{\phi_{k-1}}}{\partial{y_k}}&1  \\
 \end{array} } \right]\left[{\begin{array}{c}Y\\X\\\end{array}}\right]\\ \label{positive}
 &=&Y^TG_{k-1}(y_{k-1})Y+\left(\frac{\partial{\phi_{k-1}}}{\partial{y_k}}Y-X\right)^2. 
\end{eqnarray} 
%\end{small}
If $Y\in\R^{k-1}$ is the zero vector, $X$ must be nonzero implying that the equation (\ref{positive}) is equal to $X^2$ which is positive. On the other hand, if $Y\in\R^{k-1}$ is nonzero, $Y^TG_{k-1}(y_{k-1})Y$ is a positive scalar because $G_{k-1}$ is a Riemannian metric. Hence, $G_k$ is positive definite. Using the inequality (\ref{contsub1}), the long algebraic manipulations in (\ref{proof1}) show that $G_k$ satisfies (\ref{contractionGAS}) with the contraction rate $\lambda$. Hence, the metric $G_k$ is a contraction metric, with respect to states, for the system $\Sigma_k$.
\begin{figure*}
----------------------------------------------------------------------------------------------------------------------
\begin{small}
\begin{eqnarray}
\label{proof1}
&\left[Y^T~X\right]\left(\left(\frac{\partial\left[F_k^T+B_k^Tx_k~k_k(x,\widehat{u})\right]^T}{\partial{z_k}}\right)^TG_k(y_k)+G_k(y_k)\frac{\partial\left[F_k^T+B_k^Tx_k~k_k(x,\widehat{u})\right]^T}{\partial{z_k}}+\dot{G}_k(y_k)\right)\left[{\begin{array}{c}Y\\X\\\end{array}}\right]=\\ \notag
&\left[Y^T~X\right]\Bigg(\left[{\begin{array}{cc}
 \frac{\partial{F_k}}{\partial{y_k}}&B_k  \\
 \left(F_k+B_kx_k\right)^T\frac{\partial^2\phi_{k-1}}{\partial{y_k}^2}+\frac{\partial{\phi_{k-1}}}{\partial{y_k}}\frac{\partial{F_k}}{\partial{y_k}}+\frac{\lambda}{2}\frac{\partial{\phi_{k-1}}}{\partial{y_k}}-B_k^TG_{k-1}(y_{k-1})&-\frac{\lambda}{2}+\frac{\partial{\phi_{k-1}}}{\partial{y_k}}B_k  \\
 \end{array} }\right]^T\cdot\\ \notag&\left[\!{\begin{array}{cc}
 G_{k-1}(y_{k-1})+\left(\frac{\partial{\phi_{k-1}}}{\partial{y_k}}\right)^T\frac{\partial{\phi_{k-1}}}{\partial{y_k}}&-\left(\frac{\partial{\phi_{k-1}}}{\partial{y_k}}\right)^T  \\
 -\frac{\partial{\phi_{k-1}}}{\partial{y_k}}&1  \\
 \end{array} }\right]+\left[{\begin{array}{cc}
G_{k-1}(y_{k-1})+\left(\frac{\partial{\phi_{k-1}}}{\partial{y_k}}\right)^T\frac{\partial{\phi_{k-1}}}{\partial{y_k}}&-\left(\frac{\partial{\phi_{k-1}}}{\partial{y_k}}\right)^T  \\
-\frac{\partial{\phi_{k-1}}}{\partial{y_k}}&1  \\
\end{array} }\right]\cdot\\ \notag &\left[{\begin{array}{cc}
\frac{\partial{F_k}}{\partial{y_k}}&B_k  \\
\left(F_k+B_kx_k\right)^T\frac{\partial^2\phi_{k-1}}{\partial{y_k}^2}+\frac{\partial{\phi_{k-1}}}{\partial{y_k}}\frac{\partial{F_k}}{\partial{y_k}}+\frac{\lambda}{2}\frac{\partial{\phi_{k-1}}}{\partial{y_k}}-B_k^TG_{k-1}(y_{k-1})&-\frac{\lambda}{2}+\frac{\partial{\phi_{k-1}}}{\partial{y_k}}B_k  \\
\end{array} }\right]+ \\ \notag
&\left[\!{\begin{array}{cc}
\frac{\partial{G}_{k-1}}{\partial{y_k}}\left(F_k+B_kx_k\right)+\frac{\partial^2\phi_{k-1}}{\partial{y_k}^2}(F_k+B_kx_k)\frac{\partial\phi_{k-1}}{\partial{y_k}}+\left(\frac{\partial{\phi_{k-1}}}{\partial{y_k}}\right)^T(F_k+B_kx_k)^T\frac{\partial^2\phi_{k-1}}{\partial{y_k}^2}&-\frac{\partial^2\phi_{k-1}}{\partial{y_k}^2}(F_k+B_kx_k)  \\
-(F_k+B_kx_k)^T\frac{\partial^2\phi_{k-1}}{\partial{y_k}^2}&0 \\
\end{array} }\right]\Bigg)\cdot\\ \notag
&\left[{\begin{array}{c}Y\\X\\\end{array}}\right]=\left[Y^T~X\right]\cdot\\ \notag
&\left[\!{\begin{array}{cc}
\left(\left(\frac{\partial(F_k+B_k\phi_{k-1})}{\partial{y_k}}\right)^TG_{k-1}(y_{k-1})+G_{k-1}(y_{k-1})\frac{\partial(F_k+B_k\phi_{k-1})}{\partial{y_k}}+\frac{\partial{G}_{k-1}}{\partial{y_k}}F_k\right)-\lambda\left(\frac{\partial{\phi_{k-1}}}{\partial{y_k}}\right)^T\frac{\partial{\phi_{k-1}}}{\partial{y_k}}&\lambda\left(\frac{\partial{\phi_{k-1}}}{\partial{y_k}}\right)^T  \\
\lambda\frac{\partial{\phi_{k-1}}}{\partial{y_k}}&-\lambda  \\
\end{array} }\right]\cdot\\ \notag
&\left[{\begin{array}{c}Y\\X\\\end{array}}\right]\leq-\lambda\left[Y^T~X\right]G_k(y_k)\left[{\begin{array}{c}Y\\X\\\end{array}}\right].
\end{eqnarray}
\end{small}
----------------------------------------------------------------------------------------------------------------------
\end{figure*}
Therefore, for any $l\leq n$, the metric $G_l$ is a contraction metric, with respect to states, for the system (\ref{system3}) and with the contraction rate $\lambda$.

The proposed control law (\ref{inputor}), transforms a control system of the form (\ref{system1}) into: 
\begin{equation}
    \Sigma_n:\left\{
                \begin{array}{ll}
                 \dot{\eta}_n=F_n(\eta_n)+B_n\xi_n,\\
                 \dot{\xi}_n=k_n(\eta_n,\xi_n)+\widehat\upsilon.
                \end{array}
                \right.
\label{system2}
\end{equation}
It can be easily checked that $\widehat\upsilon$ does not appear in the variation of $\Sigma_n$ when we only have state variations. Since the metric $G_n$ is not a function of the $n$-th state, its derivative with respect to time does not include $\widehat\upsilon$. Hence, we can apply the induction results to $\Sigma_n$ to conclude that the metric $G_n$ is a contraction metric, with respect to states, for $\Sigma_n$ and with the contraction rate $\lambda$. Moreover, it can be readily seen that $G_n=\psi^*I_n$, where
\begin{equation}
\psi(x)=\left[{\begin{array}{c}x_1\\x_2-\phi_1(x)\\x_3-\phi_2(x)\\\vdots\\x_n-\phi_{n-1}(x)\\\end{array}}\right].
\label{isometry}
\end{equation}
Note that $\mathbf{d}_{I_n}$ is just the Euclidean metric and we know that $(\R^n,\mathbf{d}_{I_n})$ is a complete metric space. Moreover, since $\psi:\R^n\rightarrow\R^n$ is an isometry\footnote{Suppose $M$ and $\widetilde{M}$ are Riemannian manifolds with Riemannian metrics $G$ and $\widetilde{G}$, respectively. A smooth map
$\psi:M\rightarrow{\widetilde{M}}$ is called an isometry if it is a diffeomorphism satisfying $G=\psi^*\widetilde{G}$.}, $(\R^n,\mathbf{d}_{G_n})$ is also a complete metric space \cite{lee}. By using Theorem \ref{theoremGAS}, we conclude that a control system of the form (\ref{system1}), equipped with the state feedback control law (\ref{inputor}), is $\delta_\exists$-GAS. The $\delta_\exists$-GAS condition (\ref{delta_GAS}), as shown in \cite{aghannan}, is given by: 
\begin{equation}
\mathbf{d}_{G_n}\left(\xi_{x\upsilon}(t),\xi_{x'{\upsilon}}(t)\right)\leq e^{-\frac{\lambda}{2}{t}}\mathbf{d}_{G_n}(x,x').
\end{equation}
%where $\mathbf{d}_{G_n}(x,x')$ is the geodesic distance\footnote{The geodesic distance is the length of the shortest path between any two points on a manifold.} with respect to the metric $G_n$ between the points $x$ and $x'$. More details about the geodesic distances can be found in \cite{jost}. 
\end{proof}
\begin{remark}
The contraction metric $G_n(y_{n})$, with respect to states, for the control system (\ref{system1}), equipped with the state feedback control law (\ref{inputor}), is given by:
\begin{footnotesize}
\begin{eqnarray}
 \label{GASmetrictot}\\
 \nonumber
\left[ {\begin{array}{cc}
\left[{\begin{array}{c}
 \left[{\begin{array}{cc}
 \left[{\begin{array}{cc}
 [1]+\left(\frac{\partial\phi_1}{\partial{y}_2}\right)^T\frac{\partial\phi_1}{\partial{y}_2}&-\left(\frac{\partial\phi_1}{\partial{y}_2}\right)^T\\
 -\frac{\partial\phi_1}{\partial{y}_2}&1\\
 \end{array}}\right]+\left(\frac{\partial\phi_2}{\partial{y}_3}\right)^T\frac{\partial\phi_2}{\partial{y}_3}&-\left(\frac{\partial\phi_2}{\partial{y}_3}\right)^T\\
-\frac{\partial\phi_2}{\partial{y}_3}&1\\
 \end{array}}\right]+\cdots\\
 \vdots\\
 \end{array}}\right]+\left(\frac{\partial{\phi_{n-1}}}{\partial{y_{n}}}\right)^T\frac{\partial{\phi_{n-1}}}{\partial{y_{n}}}&-\left(\frac{\partial{\phi_{n-1}}}{\partial{y_{n}}}\right)^T  \\
 -\frac{\partial{\phi_{n-1}}}{\partial{y_{n}}}&1  \\
 \end{array} } \right],
\end{eqnarray}
\end{footnotesize}
where $y_l=[x_1,\cdots,x_{l-1}]^T$, for $l=2,\cdots,n$, and the contraction rate is $\lambda$.
\end{remark}
\begin{remark}
It can be checked that the function
\begin{equation}
\nonumber
V(x)=\frac{1}{2}\sum_{l=0}^{n-1}\left(x_{l+1}-\phi_{l}(x)\right)^2, 
\end{equation} 
is a Lyapunov function \cite{khalil} for the control system (\ref{system1}), equipped with the state feedback control law (\ref{inputor}) when $\widehat{u}=0$. Moreover, the Hessian of $V(x)$ is equal to the contraction metric $G_n$, with respect to states, defined in (\ref{GASmetrictot}).
\end{remark}

In the next theorem, we show that control law (\ref{inputor}) also enforces $\delta_\exists$-ISS.
\begin{theorem}
\label{theorem2}
For any control system $\Sigma=(\R^n,\mathsf{U},\mathcal{U},f)$ with $f$ of the form~(\ref{system1}) and for any $\lambda\in\mathbb{R}^+$, the state feedback control law:
\begin{eqnarray}
k(x,\widehat{u})&=&\frac{1}{g(x)}\bigg[k_{n}(x)-h_n(x)\bigg]+\frac{1}{g(x)}\widehat{u}, \label{inputor1}
\end{eqnarray}
where 
\begin{eqnarray}
\nonumber
k_{l}(x,\widehat{u})&=&-b_{l-1}\left(x_{l-1}-\phi_{l-2}(x)\right)-\frac{\lambda}{2}\left(x_{l}-\phi_{l-1}(x)\right)+\frac{\partial{\phi}_{l-1}}{\partial x}f(x,k(x,\widehat{u})),~\text{for}~l=1,\cdots,n,\\ \notag
\phi_{l}(x)&=&\frac{1}{b_{l}}\bigg[k_{l}(x)-h_{l}(x)\bigg],~\text{for}~l=1,\cdots,n-1,\\ \notag &&\phi_{-1}(x)=\phi_0(x)=0~~\forall x\in\R^n,~b_0=0,~\text{and}~ x_{0}=0,  \label{input1}
\end{eqnarray}
renders the control system $\Sigma$ $\delta_\exists$-ISS with respect to the input $\widehat{\upsilon}$.
\end{theorem}

\begin{proof}
Consider the following system:
%\begin{eqnarray}
%    \Sigma_{l}:\left\{
%                \begin{array}{ccl}
%                 \dot{\xi}_1&=&f_1(\xi_1)+b_1\xi_2,\\
%                 &\vdots&\\
%                 \dot{\xi}_{l-2}&=&f_{l-2}(\xi_1,\cdots,\xi_{l-2})+b_{l-2}\xi_{l-1},\\
%                 \dot{\xi}_{l-1}&=&k_{l-1}(\xi_1,\cdots,\xi_{l-1}),\\
%                \end{array}
%                \right.
%\label{system4}
%\end{eqnarray}
%or equivalently:
\begin{equation}
    \Sigma_{l}:
                \begin{array}{ll}
                 \dot{\eta}_{l}=F_{l}(\eta_{l})+B_{l}\phi_{l-1}(\eta_l),~\text{for}~l=2,\cdots,n,\\
                \end{array}
\label{system5}
\end{equation}
where \mbox{$\eta_{l}=[\xi_1,\cdots,\xi_{l-1}]^T$}, $B_{l}=[0,\cdots,0,b_{l-1}]^T\in\mathbb{R}^{l-1}$, \mbox{$y_{l}^T=[x_1,\cdots,x_{l-1}]^T$} is the state of $\Sigma_l$, and $F_{l}(\eta_{l})=[f_1(\xi,\upsilon),\cdots,f_{l-2}(\xi,\upsilon),h_{l-1}(\xi_1,\cdots,\xi_{l-1})]^T$. As proved in Theorem \ref{theorem1}, the metric $G_{l}$ in (\ref{contraction_metric}) is a contraction metric, with respect to states, for the system $\Sigma_{l}$ and with the contraction rate $\lambda$.
The proposed control law (\ref{inputor1}), transforms a control system of the form (\ref{system1}) into:
\begin{equation}
    \Sigma:\left\{
                \begin{array}{ll}
                 \dot{\eta}_{n}=F_{n}(\eta_{n})+B_{n}\xi_{n},\\
                 \dot{\xi}_n=k_n(\eta_{n},\xi_n)+\widehat{\upsilon}.
                \end{array}
                \right.
\label{systemred}
\end{equation}
Now, we show that: 
%\begin{small}
\begin{equation}
G_n(y_{n})=\left[ {\begin{array}{cc}
 G_{{n-1}}(y_{n-1})+\left(\frac{\partial{\phi_{n-1}}}{\partial{y_n}}\right)^T\frac{\partial{\phi_{n-1}}}{\partial{y_n}}&-\left(\frac{\partial{\phi_{n-1}}}{\partial{y_n}}\right)^T  \\
 -\frac{\partial{\phi_{n-1}}}{\partial{y_n}}&1  \\
 \end{array} } \right], \label{ISSmetric}
\end{equation}
%\end{small}
is a contraction metric, with respect to states and inputs, for the control system (\ref{systemred}). For $n=1$, it can be easily checked that $G_1(y_1)=1$ is a contraction metric, with respect to states and inputs, with the contraction rate $\lambda$, satisfying (\ref{contraction2ISS}) with $\alpha=2$ for the scalar control system: 
\begin{eqnarray}
\nonumber
\Sigma:\dot{\xi}_1=k_1(\xi_1)+\widehat{\upsilon}=-\frac{\lambda}{2}\xi_1+\widehat{\upsilon}. \nonumber
\end{eqnarray}
As proved in Theorem \ref{theorem1}, $G_n(y_n)$ is positive definite. Using the inequality (\ref{contsub1}) for $k=n$, long algebraic manipulations similar to those in (\ref{proof1}) show that $G_n$ satisfies (\ref{contraction2ISS}) with the contraction rate $\lambda$ and \mbox{$\alpha=2$}. Hence, the metric $G_n$ is a contraction metric, with respect to states and inputs, for the control system (\ref{systemred}).
As explained in the proof of Theorem \ref{theorem1}, we know that $(\R^n,\mathbf{d}_{G_n})$ is a complete metric space. By using Theorem \ref{theoremISS}, we conclude that a control system of the form (\ref{system1}), equipped with the state feedback control law (\ref{inputor1}), is $\delta_\exists$-ISS with respect to $\widehat{\upsilon}$. The $\delta_\exists$-ISS condition (\ref{delta_ISS}), as shown in \cite{majid}, is given by: 
\begin{eqnarray}
\notag
\mathbf{d}_{G_n}\left(\xi_{x\widehat\upsilon}(t),\xi_{x'{\widehat\upsilon'}}(t)\right)&\leq&e^{-\frac{\lambda}{2}{t}}\mathbf{d}_{G_n}(x,x')+\frac{2}{\lambda}\left(1-e^{-\frac{\lambda}{2}{t}}\right)\Vert\widehat\upsilon-\widehat\upsilon'\Vert_{\infty}\\ \notag
&\leq& e^{-\frac{\lambda}{2}{t}}\mathbf{d}_{G_n}(x,x')+\frac{2}{\lambda}\Vert\widehat\upsilon-\widehat\upsilon'\Vert_{\infty}. \label{ISScondition}
\end{eqnarray}
%where $d_{G_n}(x,x')$ is the geodesic distance with respect to the metric $G_n$ between the points $x$ and $x'$.
\end{proof}
\begin{remark}
The contraction metric, with respect to states and inputs, for a control system of the form (\ref{system1}), equipped with the state feedback control law (\ref{inputor1}), is given by (\ref{GASmetrictot}).
\end{remark}
\begin{remark}
It can be shown that the function
\begin{equation}
\nonumber
V(x)=\frac{1}{2}\sum_{l=0}^{n-1}\left(x_{l+1}-\phi_{l}(x)\right)^2, 
\end{equation} 
is an Input-to-State Stability Lyapunov function \cite{khalil} with respect to $\widehat{\upsilon}$ for a control system of the form (\ref{system1}), equipped with the state feedback control law (\ref{inputor1}). 
\end{remark}

Now, we extend the results in Theorems \ref{theorem1} and \ref{theorem2} to the class of control systems \mbox{$\Sigma=(\R^n,\R,\mathcal{U},f)$} with $f$ of the strict-feedback form \cite{miroslav}: 
\begin{eqnarray}
\label{system6}
                \begin{array}{ccl}
                 f_1(x,u)&=&h_1(x_1)+g_1(x_1)x_2,\\
                 f_2(x,u)&=&h_2(x_1,x_2)+g_2(x_1,x_2)x_3,\\
                 &\vdots&\\
                 f_{n-1}(x,u)&=&h_{n-1}(x_1,\cdots,x_{n-1})+g_{n-1}(x_1,\cdots,x_{n-1})x_n,\\
                 f_{n}(x,u)&=&h_{n}(x)+g_n(x)u,\\
                \end{array}
\end{eqnarray}
where \mbox{$x\in\mathbb{R}^{n}$} is the state and $u\in\mathsf{U}\subset\mathbb{R}$ is the control input. The functions \mbox{$h_i:\mathbb{R}^{i}\rightarrow\mathbb{R}$}, and \mbox{$g_i:\mathbb{R}^{i}\rightarrow\mathbb{R}$}, for \mbox{$i=1,\ldots,n$}, are smooth, and \mbox{$g_i(x_1,\cdots,x_i)\neq{0}$} over the domain of interest. 

In order to extend Theorems \ref{theorem1}, and \ref{theorem2} to control systems of the form (\ref{system6}), we need the following technical lemmas.

\begin{lemma}
\label{lemma1}
Let \mbox{$\Sigma=(\R^n,\mathsf{U},\mathcal{U},f)$} be a control system and let $\phi:\R^n\to\R^n$ be a smooth map with a smooth inverse. If the metric $G$ is a contraction metric, with respect to states, for $\Sigma'=(\R^n,\mathsf{U},\mathcal{U},\phi_* f)$ and with contraction rate $\lambda\in\R^+$, then the metric $\phi^*G$ is a contraction metric, with respect to states, for the system $\Sigma$ and with the contraction rate $\lambda$.
\end{lemma}

\begin{proof}
Since $G$ is a contraction metric, with respect to states, for the system ${\Sigma}'$ and with the contraction rate $\lambda$, using the inequality (\ref{dervar1}), we have:
\begin{equation}
\frac{d}{dt}\langle\delta{\eta},\delta{\eta}\rangle_{G}\leq-\lambda\langle\delta{\eta},\delta{\eta}\rangle_{G}, \label{conttran}
\end{equation}
where $\delta\eta$ is variation of the state trajectory of $\Sigma'$. Since $G$ is a metric and $\Theta(x)=\frac{\partial{\phi}}{\partial{x}}(x)$ is an invertible matrix\footnote{For any smooth map $\phi:\mathbb{R}^n\rightarrow\mathbb{R}^n$ with a smooth inverse, it is easy to show that $\frac{\partial{\phi}}{\partial{x}}(x)$ is an invertible matrix for any $x\in\mathbb{R}^n$.}, it is readily seen that $(\phi^*G)(x)$ is a positive definite matrix. We now show that the metric $\phi^*G$ is a contraction metric, with respect to states, for the system $\Sigma$.
For the coordinate transformation $\eta=\phi(\xi)$, we have:
\begin{equation}
\delta\eta=\Theta(\xi)\delta\xi. \label{vartran}
\end{equation} 
By taking the derivative of (\ref{vartran}) with respect to time, we obtain:
\begin{equation}
\frac{d}{dt}\delta\eta=\dot{\Theta}(\xi)\delta\xi+\Theta(\xi)\frac{d}{dt}\delta\xi. \label{dertran}
\end{equation} 
Using (\ref{conttran}), (\ref{vartran}), and (\ref{dertran}), we obtain:
\begin{eqnarray}
\nonumber
\frac{d}{dt}\langle\delta{\eta},\delta{\eta}\rangle_{G}&=&\left(\frac{d}{dt}\delta\eta\right)^TG\delta\eta+\delta\eta^TG\frac{d}{dt}\delta\eta+\delta\eta^T\dot{G}\delta\eta\\ \notag
&=&\left(\dot{\Theta}\delta\xi+\Theta\frac{d}{dt}\delta\xi\right)^TG\Theta\delta\xi+\left(\Theta\delta\xi\right)^TG\left(\dot{\Theta}\delta\xi+\Theta\frac{d}{dt}\delta\xi\right)+\delta\eta^T\dot{G}\delta\eta\\ \notag
&=&\left(\frac{d}{dt}\delta{\xi}\right)^T\phi^*G\delta\xi+\delta\xi^T\phi^*G\frac{d}{dt}\delta{\xi}+\delta\xi\frac{d}{dt}\left(\phi^*G\right)\delta\xi=\frac{d}{dt}\langle\delta{\xi},\delta{\xi}\rangle_{\phi^*G}\\ \notag
&\leq& -\lambda\langle\delta{\eta},\delta{\eta}\rangle_{G}=-\lambda\delta\eta^TG\delta\eta=-\lambda\left(\Theta\delta\xi\right)^TG\Theta\delta\xi=-\lambda\delta\xi^T\phi^*G\delta\xi\\ \notag
&=&-\lambda\langle\delta{\xi},\delta{\xi}\rangle_{\phi^*G}.
\end{eqnarray}
Hence, the metric $\phi^*G$ is a contraction metric, with respect to states, for the system $\Sigma$ and with the contraction rate $\lambda$.
\end{proof}
In the next lemma, we extend the results of Lemma \ref{lemma1} to contraction with respect to states and inputs.
\begin{lemma}
\label{lemma2}
Let $\Sigma=(\mathbb{R}^{n},\mathsf{U},\mathcal{U},f)$ be a control system and let $\phi:\R^n\to\R^n$ be a smooth map with a smooth inverse. If the metric $G$ is a contraction metric, with respect to states and inputs, satisfying (\ref{contraction2ISS}) with contraction rate $\lambda\in\R^+$, and $\alpha\in\R_0^+$ for $\Sigma'=(\mathbb{R}^{n},\mathsf{U},\mathcal{U},\phi_* f)$, then the metric $\phi^*G$ is a contraction metric, with respect to states and inputs, satisfying (\ref{contraction2ISS}) with the contraction rate $\lambda$, and the nonnegative constant $\alpha$ for the system $\Sigma$.
\end{lemma}

\begin{proof}
Since $G$ is a contraction metric, with respect to states and inputs, for the control system ${\Sigma}'$, satisfying (\ref{contraction2ISS}) with the contraction rate $\lambda$, and $\alpha\in\R_0^+$, using the inequality (\ref{dervar2}), we have:
\begin{equation}
\frac{d}{dt}\langle\delta{\eta},\delta{\eta}\rangle_{G}\leq-\lambda\langle\delta{\eta},\delta{\eta}\rangle_{G}+\alpha{\langle\delta{\eta},\delta{\eta}\rangle^{\frac{1}{2}}_G}\langle{\delta{\upsilon},\delta{\upsilon}}\rangle^{\frac{1}{2}}_{I_m}. \label{conttran1}
\end{equation}
Using (\ref{vartran}), (\ref{dertran}), (\ref{conttran1}), and the results of Lemma \ref{lemma1}, we obtain:
\begin{eqnarray}
\notag
\frac{d}{dt}\langle\delta{\eta},\delta{\eta}\rangle_{G}&=&\frac{d}{dt}\langle\delta{\xi},\delta{\xi}\rangle_{\phi^*G}\\ \notag
%&=&\left(\dot{\Theta}\delta\xi+\Theta\frac{d}{dt}\delta\xi\right)^TG\Theta\delta\xi+\left(\Theta\delta\xi\right)^TG\left(\dot{\Theta}\delta\xi+\Theta\frac{d}{dt}\delta\xi\right)+\delta\eta^T\dot{G}\delta\eta\\ \notag
%&=&\delta\dot{\xi}^T\left(\Theta^TG\Theta\right)\delta\xi+\delta\xi^T\left(\Theta^TG\Theta\right)\delta\dot{\xi}+\delta\xi\left(\dot{\Theta}^TG\Theta+\Theta^T\dot{G}\Theta+\Theta^TG\dot{\Theta}\right)\delta\xi \\ \notag
%&=&\frac{d}{dt}\langle\delta{\xi},\delta{\xi}\rangle_{\Theta^TG\Theta}\\ \notag
&\leq& -\lambda\langle\delta{\eta},\delta{\eta}\rangle_{G}+\alpha{\langle\delta{\eta},\delta{\eta}\rangle^{\frac{1}{2}}_G}\langle{\delta{\upsilon},\delta{\upsilon}}\rangle^{\frac{1}{2}}_{I_m}\\ \notag
&=&-\lambda\langle\delta{\xi},\delta{\xi}\rangle_{\phi^*G}+\alpha\left(\left(\Theta\delta\xi\right)^TG\Theta\delta\xi\right)^\frac{1}{2}\langle{\delta{\upsilon},\delta{\upsilon}}\rangle^{\frac{1}{2}}_{I_m}\\ \notag
&=& -\lambda\langle\delta{\xi},\delta{\xi}\rangle_{\phi^*G}+\alpha\langle\delta{\xi},\delta{\xi}\rangle^\frac{1}{2}_{\phi^*G}\langle{\delta{\upsilon},\delta{\upsilon}}\rangle^{\frac{1}{2}}_{I_m}.
\end{eqnarray}
Hence, the metric $\phi^*G$ is a contraction metric, with respect to states and inputs, satisfying (\ref{contraction2ISS}) with the contraction rate $\lambda$, and the nonnegative constant $\alpha$, for $\Sigma$. 
\end{proof}
We can now state the main result for a control system $\Sigma=(\mathbb{R}^{n},\mathsf{U},\mathcal{U},f)$ with $f$ of the form (\ref{system6}).
\begin{theorem}
\label{theorem3}
Let $\Sigma=(\mathbb{R}^{n},\mathsf{U},\mathcal{U},f)$ be a control system where $f$ is of the form~(\ref{system6}). The state feedback control law $u=k(\phi(x),\widehat{u})$, where $k$ was defined in (\ref{inputor}) and $\phi:\mathbb{R}^n\rightarrow\mathbb{R}^n$ is the smooth map (with smooth inverse) defined by:
\begin{eqnarray}
\phi(x)=\left[{\begin{array}{c}x_1\\g_1(x_1)x_2\\g_1(x_1)g_2(x_1,x_2)x_3\\\vdots\\\prod_{i=1}^{n-1}g_{i}(x_1,\cdots,x_{i})x_n\\\end{array}}\right], \label{trancor}
\end{eqnarray}
renders control system $\Sigma$ $\delta_\exists$-GAS.
\end{theorem}

\begin{proof}
The coordinate transformation \mbox{$\eta=\phi(\xi)$} transforms the control system $\Sigma=(\mathbb{R}^{n},\mathsf{U},\mathcal{U},f)$ with $f$ of the form (\ref{system6}) to the control system \mbox{$\Sigma'=(\mathbb{R}^{n},\mathsf{U},\mathcal{U},f')$}, where $f'=\phi_*f$. It can be easily checked that $f'$ has the following form:
\begin{eqnarray}
                \begin{array}{ccl}
                 f'_1(y,u)&=&h'_1(y_1)+y_2,\\
                 f'_2(y,u)&=&h'_2(y_1,y_2)+y_3,\\
                 &\vdots&\\
                 f'_{n-1}(y,u)&=&h'_{n-1}(y_1,\cdots,y_{n-1})+y_n,\\
                 f'_{n}(y,u)&=&h'_{n}(y)+g'(y)u,\\
                \end{array}
\label{system7}
\end{eqnarray}
where $h'_i:\mathbb{R}^i\rightarrow\mathbb{R}$, for $i=1,\cdots,n$, are smooth functions, $g'=\prod_{i=1}^{i=n}g_i$, and $y\in\R^n$ is the state of $\Sigma'$.
As proved in Theorem \ref{theorem1}, the state feedback control law $k$, defined in (\ref{inputor}), makes the metric $G_n$, defined in (\ref{GASmetrictot}), a contraction metric, with respect to states, for the control system $\Sigma'$ and with the contraction rate $\lambda$. 
%Consider the invertible matrix:
%%\begin{small}
%\begin{eqnarray}
%&&\Theta(x)=\frac{\partial{\Phi}}{\partial{x}}(x)=\left[\begin{array}{cccccc}1&0&0&\cdots&0\\\frac{\partial{g_1}}{\partial{x_1}}x_2&g_1(\cdot)&0&\cdots&0\\\frac{\partial{g_1g_2}}{\partial{x_1}}x_3&\frac{\partial{g_1g_2}}{\partial{x_2}}x_3&g_1(\cdot)g_2(\cdot)&\cdots&0\\\vdots&\vdots&&\ddots&\vdots\\\frac{\partial{\prod_{i=1}^{n-1}g_{i}}}{\partial{x_1}}x_n&\frac{\partial{\prod_{i=1}^{n-1}g_{i}}}{\partial{x_2}}x_n&\cdots&&\prod_{i=1}^{n-1}g_{i}(\cdot)\\\end{array}\right].
%\end{eqnarray}
%\end{small}
%is invertible because all the elements on the diagonal are nonzero over the domain of the interest. 
As proved in Lemma \ref{lemma1}, the metric $\phi^*G_n$ is a contraction metric, with respect to states, with the contraction rate $\lambda$, for the control system $\Sigma$, equipped with the state feedback control law $k(\phi(x),\widehat{u})$. Since $(\R^n,\mathbf{d}_{G_n})$ is a complete metric space and $\phi$ is an isometry, $(\R^n,\mathbf{d}_{\phi^*G_n})$ is also a complete metric space \cite{lee}. Therefore, the state feedback control law $k(\phi(x),\widehat{u})$ makes the control system $\Sigma$ $\delta_\exists$-GAS.
\end{proof}
The $\delta_\exists$-ISS version of Theorem \ref{theorem3} is given by the following result.
\begin{theorem}
\label{theorem4}
Let $\Sigma=(\mathbb{R}^{n},\mathsf{U},\mathcal{U},f)$ be a control system where $f$ is of the form~(\ref{system6}). The state feedback control law $u=k(\phi(x),\widehat{u})$, where $k$ and $\phi$ were defined in (\ref{inputor1}) and (\ref{trancor}), respectively, 
%$\phi:\mathbb{R}^n\rightarrow\mathbb{R}^n$ is the smooth map (with smooth inverse) defined by:
%\begin{eqnarray}
%\phi(x)=\left[{\begin{array}{c}x_1\\g_1(x_1)x_2\\g_1(x_1)g_2(x_1,x_2)x_3\\\vdots\\\prod_{i=1}^{n-1}g_{i}(x_1,\cdots,x_{i})x_n\\\end{array}}\right], \label{trancor}
%\end{eqnarray}
renders control system $\Sigma$ $\delta_\exists$-ISS with respect to the input $\widehat{\upsilon}$.
\end{theorem}

\begin{proof}
By following the same steps as in the proof of Theorem \ref{theorem3}, and using Lemma \ref{lemma2}, we obtain that the state feedback control law $u=k(\phi(x),\widehat{u})$ makes the metric $\phi^*G_n$ a contraction metric, with respect to states and inputs, for $\Sigma$ and with the contraction rate $\lambda$. Hence, the control system $\Sigma$, equipped with the state feedback control law $k(\phi(x),\widehat{u})$, is $\delta_\exists$-ISS with respect to the input $\widehat{\upsilon}$.
\end{proof}
\begin{remark}
Although we only discussed single input control systems, extensions to multi input control systems are straightforward using the techniques in \cite{miroslav}.
\end{remark}

\section{Backstepping Controller Design for a synchronous generator}
We illustrate the results in this paper on a synchronous generator \cite{roosta} connected through a transmission line to an infinite bus. 
%\begin{figure}
%\begin{center}
%\includegraphics[width=8cm]{generator}
%\caption {Generator connected through a transmission line to an infinite bus.}\label{fig1}
%\end{center}
%\end{figure}
The control system $\Sigma=\left(\R^3,\mathsf{U},\mathcal{U},f\right)$ with $f$ of the form:
\begin{eqnarray}\nonumber
                 f_1(x,u)&=&x_2,\\ \label{generator}
                 f_2(x,u)&=&-Ex_2+FP_{m0}+V_sGe_{q0}\sin(\delta_0+x_1)+V_sG\sin(\delta_0+x_1)x_3, \\ \notag
                 f_3(x,u)&=&-Ix_3+JV_s\sin(\delta_0+x_1)x_2-Ie_{q0}+IK_cu,
\end{eqnarray}
%where
%\begin{eqnarray}
%\nonumber
%h_1(x_1)&=&0,\\ \notag
%g_1(x_1)&=&1,\\ \notag
%h_2(x_1,x_2)&=&-Ex_2+FP_{m0}+V_sGe_{q0}\sin(\delta_0+x_1),\\ \notag
%g_2(x_1,x_2)&=&V_sG\sin(\delta_0+x_1),\\ \notag
%h_3(x)&=&-Ix_3+JV_s\sin(\delta_0+x_1)x_2-Ie_{q0}, \\ \notag
%g_3(x)&=&IK_c,
%\end{eqnarray}
models a synchronous generator connected to an infinite bus. In the aforementioned model, $x_1$ is the deviation of the power angle, $x_2$ is the relative speed of the rotor of the generator, $x_3$ is the deviation of the quadrature axis voltage of the generator, $\delta_0$ is the operating point of the power angle, $P_{m0}$ is the operating point of the mechanical input power, $e_{q0}$ is the operating point of the quadrature axis voltage of the generator, $V_s$ is the infinite bus voltage, $K_c$ is the gain of the excitation amplifier, and $u$ is the input of the silicon-controlled rectifier amplifier of the generator. Other  parameters in (\ref{generator}) are given by: $E=\frac{D}{2H}$, $I=\frac{1}{T'}$, $F=\frac{\omega_0}{2H}$, $G=-\frac{\omega_0}{2H}\frac{1}{X_{qs}}$, and $J=\frac{X_q-X'_d}{X'_{ds}}$, where $D$ is the per-unit damping constant, $H$ is the inertia constant, $\omega_0$ is the synchronous generator speed, $T'=\frac{X'_{ds}}{X_{qs}}T'_{d0}$, $X_{qs}=X_T+\frac{1}{2}X_L+X_q$, $X'_{ds}=X_T+\frac{1}{2}X_L+X'_d$, $T'_{d0}$ is the direct axis transient short-circuit time constant, $X_T$ is the reactance of the transformer, $X_q$ is the quadrature axis reactance, $X'_{d}$ is the direct axis transient reactance and $X_L$ is the reactance of the transmission line. We assume that $\sin(\delta_0+x_1)$ is nonzero over the domain of the interest.

The control system (\ref{generator}) is of the form (\ref{system6}). The coordinate transformation (\ref{trancor}), given by:
\begin{equation}
\left[{\begin{array}{c}\eta_1\\\eta_2\\\eta_3\\\end{array}}\right]=\phi(\xi)=\left[{\begin{array}{c}\xi_1\\\xi_2\\V_sG\sin(\delta_0+\xi_1)\xi_3\\\end{array}}\right], \label{trancorex}
\end{equation}
transforms the control system $\Sigma=\left(\R^3,\mathsf{U},\mathcal{U},f\right)$ to the control system $\Sigma'=\left(\R^3,\mathsf{U},\mathcal{U},f'\right)$ with $f'=\phi_*f$ of the form:
\begin{eqnarray}\nonumber
                  f'_1(y,u)&=&h'_1(y_1)+y_2=y_2,\\ \label{generatortrancor}
                 f'_2(y,u)&=&h'_2(y_1,y_2)+y_3=-Ey_2+FP_{m0}+V_sGe_{q0}\sin(\delta_0+y_1)+y_3, \\ \notag
                 f'_3(y,u)&=&h'_3(y)+g'(y)u=-IV_sGe_{q0}\sin(\delta_0+y_1)+JV_s^2G\sin^2(\delta_0+y_1)y_2-Iy_3+\\ \notag
                 &&\cot(\delta_0+y_1)y_2y_3+IK_cV_sG\sin(\delta_0+y_1)u.
\end{eqnarray}
%where
%\begin{eqnarray}
%\nonumber
%h'_1(y_1)&=&0,\\ \notag
%h'_2(y_1,y_2)&=&-Ey_2+FP_{m0}+V_sGe_{q0}\sin(\delta_0+y_1),\\ \notag
%h'_3(y)&=&-IV_sGe_{q0}\sin(\delta_0+y_1)+JV_s^2G\sin^2(\delta_0+y_1)y_2-Iy_3+\cot(\delta_0+y_1)y_2y_3, \\ \notag
%g'_3(y)&=&IK_cV_sG\sin(\delta_0+y_1).
%\end{eqnarray}

By using the results in Theorem \ref{theorem2} for a control system of the form (\ref{generatortrancor}) and for $\lambda=2$, we have:
\begin{eqnarray}
\nonumber
\phi_1(\eta_1)&=&-\eta_1,\\ \notag
\phi_2(\eta_1,\eta_2)&=&-2\eta_1+\left(E-2\right)\eta_2-FP_{m0}-V_sGe_{q0}\sin(\delta_0+\eta_1), \\ \notag
k_3(\eta)&=&\left(-5+3E-E^2\right)\eta_2-3\eta_1+\left(E-3\right)\eta_3+\left(E-3\right)FP_{m0}\\ \notag
&&+\left(E-3\right)V_sGe_{q0}\sin(\delta_0+\eta_1)-V_sGe_{q0}\cos(\delta_0+\eta_1)\eta_2.
\end{eqnarray}
Therefore, the state feedback control law:
\begin{eqnarray}
\label{controlinput}
k(\eta,\widehat{\upsilon})&=&\frac{1}{g'(\eta)}\left[k_3(\eta)-h'_3(\eta)\right]+\frac{1}{g'(\eta)}\widehat{\upsilon}\\ \notag
&=&\frac{1}{IK_cV_sG\sin(\delta_0+\eta_1)}\bigg[\left(-5+3E-E^2\right)\eta_2-3\eta_1+\left(E-3+I\right)\eta_3+\left(E-3\right)FP_{m0}+\\ \notag
&&\left(E-3+I\right)V_sGe_{q0}\sin(\delta_0+\eta_1)-V_sGe_{q0}\cos(\delta_0+\eta_1)\eta_2-JV_s^2G\sin^2(\delta_0+\eta_1)\eta_2-\\ \notag
&&\cot(\delta_0+\eta_1)\eta_2\eta_3\bigg]+\frac{\widehat{\upsilon}}{IK_cV_sG\sin(\delta_0+\eta_1)},
\end{eqnarray} 
makes the control system $\Sigma'$ $\delta_\exists$-ISS with respect to the input $\widehat{\upsilon}$. The corresponding contraction metric, with respect to states and inputs, for the control system (\ref{generatortrancor}) is given by:
\begin{footnotesize}
\begin{eqnarray}
\notag
&G(y)=\left[ {\begin{array}{cc}
\left[ {\begin{array}{cc}
 1+\left(\frac{\partial{\phi_1}}{\partial{y_1}}\right)^T\frac{\partial{\phi_1}}{\partial{y_1}}&-\left(\frac{\partial{\phi_1}}{\partial{y_1}}\right)^T  \\
 -\frac{\partial{\phi_1}}{\partial{y_1}}&1 \\
 \end{array} } \right]+\left(\frac{\partial{\phi_2}}{\partial{z_2}}\right)^T\frac{\partial{\phi_2}}{\partial{z_2}}&-\left(\frac{\partial{\phi_2}}{\partial{z_2}}\right)^T  \\
 -\frac{\partial{\phi_2}}{\partial{z_2}}&1 \\
 \end{array} } \right]=\\ \notag
 &\left[\begin{array}{ccc}
 2+\left(2+V_sGe_{q0}\cos(\delta_0+y_1)\right)^2&-2E+5-\left(E-2\right)V_sGe_{q0}\cos(\delta_0+y_1)&2+V_sGe_{q0}\cos(\delta_0+y_1)\\
 -2E+5-\left(E-2\right)V_sGe_{q0}\cos(\delta_0+y_1)&\left(E-2\right)^2+1&2-E\\
 2+V_sGe_{q0}\cos(\delta_0+y_1)&2-E&1\\
 \end{array}\right],\notag
\end{eqnarray}
\end{footnotesize}
where $z_2^T=[y_1~y_2]^T$.
By using Theorem \ref{theorem4}, the state feedback control law (\ref{controlinput}), and the coordinate transformation (\ref{trancorex}), we obtain the state feedback control law $k(\phi(x),\widehat{u})$ making $\Sigma$ $\delta_\exists$-ISS with respect to the input $\widehat\upsilon$. The corresponding contraction metric, with respect to states and inputs, for the control system $\Sigma$ is given by:
\begin{equation}
\left(\phi^*G\right)(x)=\Theta^T(x)G(\phi(x))\Theta(x), \notag
\end{equation}
where
\begin{eqnarray}
\Theta(x)=\frac{\partial{\phi}}{\partial{x}}(x)=\left[ {\begin{array}{ccc}
1&0&0  \\
0&1&0 \\
V_sG\cos(\delta_0+x_1)x_3&0&V_sG\sin(\delta_0+x_1)\\
 \end{array} } \right]. 
\end{eqnarray}
Since the map $\phi$ does not transform the first and the second coordinates and the metric $G$ is only a function of the first coordinate, we have $G(\phi(x))=G(x)$.

\section{Discussion}
In this paper we developed a backstepping procedure to design controllers enforcing incremental stability. Where before we could apply backstepping to construct stabilizing controllers, we can now apply the results in this paper to construct incrementally stabilizing controllers.

\bibliographystyle{alpha}
\bibliography{reference}

\begin{thebibliography}{PPvdWN04}

\bibitem[Ang02]{angeli}
D.~Angeli.
\newblock A {L}yapunov approach to incremental stability properties.
\newblock {\em IEEE Transactions on Automatic Control}, 47(3):410--21, 2002.

\bibitem[AR03]{aghannan}
N.~Aghannan and P.~Rouchon.
\newblock An intrinsic observer for a class of {L}agrangian systems.
\newblock {\em IEEE Transactions on Automatic Control}, 48(6):936--945, 2003.

\bibitem[AS99]{sontag}
D.~Angeli and E.~D. Sontag.
\newblock Forward completeness, unboundedness observability, and their
  {L}yapunov characterizations.
\newblock {\em Systems and Control Letters}, 38:209--217, 1999.

\bibitem[CS87]{crouch}
P.~E. Crouch and A.~J. Van~Der Schaft.
\newblock {\em Variational and hamiltonian control systems}.
\newblock Springer, 1987.

\bibitem[Dem61]{demidovich1}
B.~P. Demidovich.
\newblock Dissipativity of a nonlinear system of differential equations.
\newblock {\em ser. matem. mekh. (in Russian)}, 6:19--27, 1961.

\bibitem[Dem67]{demidovich}
B.~P. Demidovich.
\newblock {\em Lectures on stability theory (in Russian)}.
\newblock Nauka, Moscow, 1967.

\bibitem[GPT09]{girard2}
A.~Girard, G.~Pola, and P.~Tabuada.
\newblock Approximately bisimilar symbolic models for incrementally stable
  switched systems.
\newblock {\em IEEE Transactions on Automatic Control}, 55(1):116--126, January
  2009.

\bibitem[JL02]{jouffroy1}
J.~Jouffroy and J.~Lottin.
\newblock Integrator backstepping using contraction theory: a brief
  methodological note.
\newblock {\em in Proceedings of 15th IFAC World Congress}, 2002.

\bibitem[Jou05]{jerome}
J.~Jouffroy.
\newblock Some ancestors of contraction analysis.
\newblock {\em Proceedings of the 44th IEEE Conference on Decision and
  Control}, pages 5450--5455, December 2005.

\bibitem[Kha96]{khalil}
H.~K. Khalil.
\newblock {\em Nonlinear systems}.
\newblock Prentice-Hall, Inc., New Jersey, 2nd edition, 1996.

\bibitem[KKK95]{miroslav}
M.~Krstic, I.~Kanellakopoulos, and P.~P. Kokotovic.
\newblock {\em Nonlinear and adaptive control design}.
\newblock John Wiley and Sons, 1995.

\bibitem[Lee03]{lee}
J.~M. Lee.
\newblock {\em Introduction to Smooth Manifolds}.
\newblock Springer-Verlag, 2003.

\bibitem[LS98]{lohmiller}
W.~Lohmiller and J.~J. Slotine.
\newblock On contraction analysis for non-linear systems.
\newblock {\em Automatica}, 34(6):683--696, 1998.

\bibitem[PGT08]{pola}
G.~Pola, A.~Girard, and P.~Tabuada.
\newblock Approximately bisimilar symbolic models for nonlinear control
  systems.
\newblock {\em Automatica}, 44(10):2508--2516, 2008.

\bibitem[PPvdWN04]{pavlov1}
A.~Pavlov, A.~Pogromvsky, N.~van~de Wouw, and H.~Nijmeijer.
\newblock Convergent dynamics, a tribute to boris pavlovich demidovich.
\newblock {\em Systems and Control Letters}, 52:257--261, 2004.

\bibitem[PT09]{pola1}
G.~Pola and P.~Tabuada.
\newblock Symbolic models for nonlinear control systems: alternating
  approximate bisimulations.
\newblock {\em SIAM Journal on Control and Optimization}, 48(2):719--733,
  February 2009.

\bibitem[PvdWN05]{pavlov}
A.~Pavlov, N.~van~de Wouw, and H.~Nijmeijer.
\newblock {\em Uniform output regulation of nonlinear systems: a convergent
  dynamics approach}.
\newblock Springer, Berlin, 2005.

\bibitem[PvdWN07]{pavlov2}
A.~Pavlov, N.~van~de Wouw, and H.~Nijmeijer.
\newblock Global nonlinear output regulation: convergence-based controller
  design.
\newblock {\em Automatica}, 43:456--463, January 2007.

\bibitem[RdB09]{russo}
G.~Russo and M.~di~Bernardo.
\newblock Contraction theory and the master stability function: linking two
  approaches to study synchnorization in complex networks.
\newblock {\em IEEE Transactions on Circuit and Systems II}, 56:177--181, 2009.

\bibitem[RdBS09]{russo1}
G.~Russo, M.~di~Bernardo, and E.~D. Sontag.
\newblock Global entrainment of transcriptional systems to periodic inputs.
\newblock {\em Submitted for publication}, arXiv: 0907.0017, 2009.

\bibitem[RGHS01]{roosta}
A.~R. Roosta, D.~Georges, and N.~Hadj-Said.
\newblock Nonlinear control for power systems based on backstepping method.
\newblock {\em in Proceedings of the 40th IEEE Conference on Decision and
  Control}, pages 3037--3042, December 2001.

\bibitem[SK08a]{sharma2}
B.~B. Sharma and I.~N. Kar.
\newblock Adaptive control of wing rock system in uncertain environment using
  contraction theory.
\newblock {\em in Proceedings of IEEE American Control Conference}, pages
  2963--2968, June 2008.

\bibitem[SK08b]{sharma1}
B.~B. Sharma and I.~N. Kar.
\newblock Design of asymptotically convergent frequency estimator using
  contraction theory.
\newblock {\em IEEE Transactions on Automatic Control}, 53(8):1932--1937,
  September 2008.

\bibitem[SK09]{sharma}
B.~B. Sharma and I.~N. Kar.
\newblock Contraction based adaptive control of a class of nonlinear systems.
\newblock {\em in Proceedings of IEEE American Control Conference}, pages
  808--813, June 2009.

\bibitem[Son98]{sontag1}
E.~D. Sontag.
\newblock {\em Mathematical control theory}, volume~6.
\newblock Springer-Verlag, New York, 2nd edition, 1998.

\bibitem[vdWP08]{wouw}
N.~van~de Wouw and A.~Pavlov.
\newblock Tracking and synchronisation for a class of pwa systems.
\newblock {\em Automatica}, 44:2909--2915, 2008.

\bibitem[WS05]{wang}
W.~Wang and J.~J. Slotine.
\newblock On partial contraction analysis for coupled nonlinear oscillators.
\newblock {\em Biological Cybernetics}, 92:38--53, 2005.

\bibitem[Zam63]{zames1}
G.~Zames.
\newblock Functional analysis applied to nonlinear feedback systems.
\newblock {\em IEEE Transaction on Circuit Theory}, 10(3):392--404, September
  1963.

\bibitem[Zam96]{zames}
G.~Zames.
\newblock Input-output feedback stability and robustness, 1959-85.
\newblock {\em IEEE Control Systems Magazine}, 16(3):61--66, 1996.

\bibitem[ZPJT10]{majid}
M.~Zamani, G.~Pola, M.~Mazo Jr., and P.~Tabuada.
\newblock Symbolic models for nonlinear control systems without stability
  assumptions.
\newblock {\em Submitted for publication}, arXiv:1002.0822, February 2010.

\end{thebibliography}
\end{document}